\date{}
\documentclass[11pt]{article}
\usepackage{hyperref}
\usepackage{tikz}

\usepackage{tikz-cd}

\usepackage{amsmath,amsthm,amsfonts,amssymb,amsopn,amscd,bm} 
\usepackage{color}

\def\qed{{\unskip\nobreak\hfil\penalty50
\hskip2em\hbox{}\nobreak\hfil$\square$
\parfillskip=0pt \finalhyphendemerits=0\par}\medskip}
\def\proof{\trivlist \item[\hskip \labelsep{\bf Proof.\ }]}
\def\endproof{\null\hfill\qed\endtrivlist\noindent}

\def\tilde{\widetilde}

\def\a{\alpha}
\def\b{\beta}

\def\l{\lambda}

\def\D{{\cal D}}
\def\E{{\cal E}}
\def\F{{\cal F}}

\def\M{{\cal M}}

\def\R{{\cal R}}

\def\H{{\cal H}}

\def\S{{\cal S}}

\def\f{{\varphi}}
\def\s{{\sigma}}

\def\l{{\lambda}}

\def\PSL{{{\rm PSL}(2,\mathbb R)}}

\def\S2{S^{1(2)}}

\def\RR{\mathbb R}

\newtheorem{theorem}{Theorem}[section]
\newtheorem{lemma}[theorem]{Lemma}

\newtheorem{corollary}[theorem]{Corollary}

\newtheorem{proposition}[theorem]{Proposition}

\theoremstyle{remark} 

\newcommand{\ben}{\begin{equation}}
\newcommand{\een}{\end{equation}}

\newcommand{\bthm}{\begin{theorem}}
\newcommand{\ethm}{\end{theorem}}

\newcommand{\bprop}{\begin{proposition}}
\newcommand{\eprop}{\end{proposition}}
\newcommand{\bcor}{\begin{corollary}}
\newcommand{\ecor}{\end{corollary}}
\newcommand{\blem}{\begin{lemma}}
\newcommand{\elem}{\end{lemma}}

\def\PSL{PSU(1,1)}

\def\SL2{{{\rm SL}(2,\R)}}

\def\PSL2{{{\rm PSL}(2,\Reali)}}

\def\U1{{{\rm V}(1)}}
\def\SU2{{{\rm SV}(2)}}

\def\SU{{{\rm SU}}}

\def\D{{\mathcal D}}
\def\F{{\mathcal F}}
\def\H{{\mathcal H}}

\def\M{{\mathcal M}}

\def\S{{\mathcal S}}
\def\T{{\mathcal T}}

\def\bx{{\bf x}}

\def\RR{{\mathbb R}}

\def\a{\alpha}
\def\b{\beta}
        
\def\d{\delta}        
  
\def\l{\lambda}

\def\s{\sigma}

\def\f{\varphi}

\def\S{{S(\RR^d)}}
\def\Sr{{S_{\rm r}(\RR^d)}}
\def\ov{\overline}

\title{\Huge{Signal communication and modular theory}}

\author{{\sc Roberto Longo\thanks{
E-mail: {longo@mat.uniroma2.it}
}
}\\
Dipartimento di Matematica,
Universit\`a di Roma Tor Vergata,\\
Via della Ricerca Scientifica, 1, I-00133 Roma, Italy
}

\begin{document}

\maketitle

\begin{abstract}
We propose a conceptual frame to interpret the prolate differential operator, which appears in Communication Theory, as an entropy operator; indeed, we write its expectation values as a sum of terms, each subject to an entropy reading by an embedding suggested by Quantum Field Theory. This adds meaning to the classical work by Slepian et al. on the problem of simultaneously concentrating a function and its Fourier transform, in particular to the ``lucky accident" that the truncated Fourier transform commutes with the prolate operator.  The key is the notion of entropy of a vector of a complex Hilbert space with respect to a real linear subspace, recently introduced by the author by means of the Tomita-Takesaki modular theory of von Neumann algebras. We consider a generalization of the prolate operator to the higher dimensional case and show that it admits a natural extension commuting with the truncated Fourier transform; this partly generalizes the one-dimensional result by Connes to the effect that there exists a natural selfadjoint extension to the full line commuting with the truncated Fourier transform. 
\end{abstract}

\newpage

\section{Introduction}
The aim of this paper is to provide an interpretation of the prolate operator, which plays an important role in the theory of signal transmission, as an entropy operator, by means of the modular theory of von Neumann algebras, following recent concepts and abstract analysis of entropy in the framework of Quantum Field Theory. We begin with a brief account of the background of our work. 
\smallskip

\noindent
{\it Band limited signals.}
Suppose Alice sends a signal to Bob that is codified by a function of time $f$. Bob can measure the value $f$ only within a certain time interval; moreover, the frequency of $f$ is filtered by the signal device within a certain interval. For simplicity, let us assume these intervals are both equal to the interval $B=(-1,1)$. As is well known, if a function $f$ and its Fourier transform $\hat f$  are both supported in bounded intervals, then $f$ is the zero function. So one is faced with the problem of simultaneously maximizing the portions of energy and amplitude spectrum within the intervals
\ben\label{SPe}
||f||^2_B/ ||f||^2\, ,\quad ||\hat f||^2_B/ ||\hat f||^2\, ,
\een
where $||\cdot||, ||\cdot||_B$ denote the $L^2$-norms on $\RR$ and $B$, the {\it concentration problem}. 

The problem of best approximating, with support concentration, a function and its Fourier transform is a classical problem; in particular, it lies behind Heisenberg uncertainty relations in Quantum Mechanics and is studied in Quantum Field Theory too, see \cite{Ja}. 

In the `60ies, this problem was studied in seminal works by Slepian, Pollak and Landau  \cite{SP1,LP}, see also \cite{S83}. The functions that best maximize \eqref{SPe} are eigenfunctions of the angle operator associated with the truncated Fourier transform. This is a Hilbert-Schmidt integral operator whose spectral analysis is not easily doable a priori. However, by the {\it lucky accident} figured out in \cite{SP1}, this integral operator commutes with a linear differential operator, the {\it prolate operator}
\ben\label{prog}
W = \frac{d}{dx}(1 -x^2)\frac{d}{dx}   - x^2 \, ,
\een
that shares its eigenfunctions with the angle operator, so these eigenfunctions were computed. 

$W$ is a classical operator, it arises by separating the 3-dimensional scalar wave
equation in a prolate spheroidal coordinate system. More recently, Connes has reconsidered and raised new interest in this operator \cite{Co}.
The papers \cite{CMa, CMo} show an impressive relation of the prolate spectrum with the asymptotic distribution of the zeros of the Riemann $\zeta$-function. Our paper is not related to this point; however, our Sect. \ref{Prol} is inspired and generalizes a small part of the analysis in \cite{CMo}.  

Our purpose is to understand the role of the prolate operator on a conceptual basis, in relation to the mentioned lucky accident. We shall argue that the prolate operator gives rise to an {\it entropy operator}, in a sense that will be explained. 
Within our aim, we shall generalize the prolate operator in higher dimensions and analyze it guided by the Quantum Field Theory context.

We shall consider the prolate operator
\ben\label{Pro1}
W_{\rm min} = (1 -r^2)\nabla^2 - 2r\partial_r   - r^2 \, ,
\een
on the Schwarz space $\S$, with $r$ the radial coordinate in $\RR^d$, and show it admits a natural closed extension $W$ that commutes with the truncated Fourier transform. We shall see that the expectation values of $\pi E_B W$ on $L^2(\RR^d)$, with $E_B$ the orthogonal projection onto $L^2(B)$, is positive, selfadjoint and its expectation values are indeed entropy quantities. 

In the one-dimensional case, $W$ itself is selfadjoint \cite{CMo}, and this is probably true also in higher dimensions; however, for our aim, it suffices to know that $E_B W$ is selfadjoint. 
\smallskip

\noindent
{\it Modular theory, the entropy of a vector.}
In the `70ies, Tomita and Takesaki uncovered a fundamental, deep operator algebraic structure. In particular, associated with any faithful, normal state $\f$ of a von Neumann algebra $\M$, there is a canonical one-parameter automorphism group $\s^\f$ of $\M$, the {\it modular group}, see \cite{Tak}. The relevance of this intrinsic evolution in Physics was soon realized in the framework of Quantum Statistical Mechanics since $\s^\f$ is characterized by the KMS thermal equilibrium condition, see \cite{H}. 	

Now, part of the modular theory shows up at a more elementary level, with potential points of contact with contexts not immediately related to Operator Algebras: the general framework is simply provided by a real linear subspace of complex Hilbert space, cf. \cite{LRT,RvD}. 

Let $\H$ be a complex Hilbert space and $H$ a real linear subspace of $\H$; by considering its closure, we may assume that $H$ is closed.
$H$ is said to be a {\it standard subspace} if $H$ is closed and
$\overline{H + iH} = \H$ , $H\cap i H = \{0\} $. 
Every closed real linear subspace $H$ has a standard subspace direct sum component and we may assume that $H$ is standard by restricting to this component. 

With $H$ standard, the anti-linear operator $S: H + iH \to H + iH$, $S(\Phi_1 + i\Phi_2) = \Phi_1 - i\Phi_2$ is then  well-defined, closed, involutive. Its polar decomposition $S = J_H\Delta_H^{1/2}$ then gives an anti-linear, involutive unitary $J_H$ and a positive, non-singular, selfadjoint operator $\Delta_H$ on $\H$, the {\it modular conjugation} and the {\it modular operator}, such that
\[
\Delta_H^{is} H = H\ , \quad J_H H = H' \ ,
\]
$s\in\RR$;
here $H'$ is the symplectic complement $H' = (iH)^{\bot_\mathbb R}$ of $H$, the orthogonal of $iH$ with respect to the real scalar product $\Re(\cdot,\cdot)$.
We refer to  \cite{LN} for the modular theory and basic results on standard subspaces. 
 
We say that the standard subspace $H$ is {\it factorial} if
$H\cap H' =\{0\}$. 
Thus $H +H'$ is dense in $\H$ and $H+ H'$ is the direct sum (as linear space) of $H$ and $H'$. 
Again, we may assume that $H$ is factorial by restricting to the factorial component. 
Our abstract results have an immediate extension to the non-factorial, non-standard case. 

The {\it cutting projection} relative to $H$ is the real linear, densely defined projection 
\[
P_H : H + H' \to H\, , \ \Phi + \Phi' \mapsto \Phi\, .
\]
The {\it entropy of a vector} $\Phi\in\H$ with respect to a standard subspace $H\subset\H$ is defined by
\ben\label{Se}
S_\Phi = \Im(\Phi, P_H i\log\Delta_H\, \Phi) = (\Phi, i P_H i\log\Delta_H\, \Phi) \, ;
\een
this notion was introduced in \cite{L19,CLR}. 
A first way to realize the entropy meaning of $S_\Phi$ is to consider the von Neumann algebra $R(H)$ associated with $H$ by the second quantization on the Fock Hilbert space over $\H$; then $S_\Phi$ is Araki's relative entropy  \cite{Ar76}  between the coherent state associated with $\Phi$ and the vacuum state on $R(H)$. However, in this paper, this fact does not play any direct role. 

Note that $ i P_H i\log\Delta_H$ is a real linear operator. This is our first instance of an {\it entropy operator}, namely a real linear, positive, selfadjoint operator whose expectation values give the entropy of states. In concrete situations, the subspace $H$ may correspond to a region of a manifold and $\Phi$ to signal, then $S_\Phi$ acquires the meaning of local entropy of $\Phi$. 
\smallskip

\noindent
{\it Entropy density of a wave packet.}
The local entropy of a wave packet has been studied in \cite{L19,CLR,CLRR} for the case of a half-space, and in \cite{LM} for the space ball case, which is directly related to the present paper; these works were motivated by Quantum Field Theory. 

Let $\T$ be the real linear space of wave packets, that is $\Phi\in\T$ if $\Phi$ is a real function on $\RR^{1+d}$ that satisfies the wave equation
$\partial_t^2 \Phi = \nabla_x^2 \Phi$, with Cauchy data in the real Schwarz space $\Sr$. Quantum Relativistic Mechanics tells us that $\T$ is equipped with a natural (Lorentz invariant) complex pre-Hilbert structure so, by completion, we get a complex Hilbert space $\H$. Wave packets with Cauchy data supported in the open, unit ball $B$ of $\RR^d$ form a real linear subspace of $\H$ denoted by $H= H(B)$ (after closure). The entropy of $\Phi$ in $B$ is given by
\ben\label{Sform}
S_{\Phi} = \pi\int_{B} (1 - r^2) \langle T_{00}\rangle_{\Phi}\, dx 
+  \pi D\int_{B} \Phi^2 dx \, .
\een
Here $D = (d-1)/2$ and  $ \langle T_{00}\rangle_{\Phi} = \frac12 \big( (\partial_t \Phi)^2  + |\nabla_x \Phi|^2    \big)$ is the energy density
of $\Phi$. We discuss here $d>1$ case; the case $d=1$ is similar but requires modifications due to infrared singularities, which is not important for our discussion. 

The two terms in $ \langle T_{00}\rangle_{\Phi} $ have separate meanings, they correspond to the kinetic and to the potential energy of the wave packet. $\H$ is naturally a direct sum of the two real Hilbert subspaces associated with the Cauchy data. 

In terms of the Cauchy data $f,g$ of $\Phi$, the modular Hamiltonian $\log\Delta_B$ relative to $B$ is given by
\ben
\imath\log\Delta_B =
\pi \left[
\begin{matrix}
0 & M \\ L - 2D& 0
\end{matrix}\right]
 = \pi \left[
\begin{matrix}
0 & (1 - r^2) \\ (1 - r^2)\nabla^2 - 2r\partial_r  - 2D & 0
\end{matrix}\right]\, ,
\een
\cite{LM}. 
Here,
\ben\label{Le}
L = (1 - r^2)\nabla^2 - 2r\partial_r 
\een
is a higher-dimensional Legendre operator. 

Each of the two terms in the expression of $S_\Phi$,
\[
S_\Phi = -\pi (f,  L_D f)_B + \pi (g,  Mg)_B \, ,
\]
$L_D\equiv L - 2D$,
have an entropy meaning. As we will discuss on general grounds, $-\pi (f,  L_D f)_B$ is the {\it field entropy} of $f$, and $\pi(g,  Mg)_B$ is the {\it momentum, or parabolic, entropy} of $g$, in $B$. We infer that also $-\pi (f,  L f)_B$ is an entropy quantity, the {\it Legendre entropy} of $f$ in $B$. 
\smallskip

\noindent
{\it The measure of concentration.}
We now return to the Communication Theory setting. The truncated Fourier transform operator is obviously defined in any space dimension. Indeed, the concentration problem often arises in higher dimensions too.  It is also studied in \cite{S4}, although with a point of view different from the one in this paper. 

As said, the higher dimensional prolate operator \eqref{Pro1} extends to a natural operator $W $ on $L^2(\RR^d)$, that commutes both with the Fourier and the truncated Fourier transforms; $W$ also commutes with the orthogonal projection $E_B$ onto $L^2(B)$ and its Fourier conjugate $\hat E_B$.  

As $-W + M = - L +1$, given $f\in \S$ real, we have
\[
- \pi (f, W f)_B + \pi (f, Mf)_B = - \pi (f, L  f)_B  +  \pi (f,f)_B\, ;
\]
that is, $- \pi (f, W f)_B$ is the sum of the Legendre entropy of $f$ and $\pi ||f||^2_B$ (that we call the Born entropy), minus the parabolic  entropy of $f$, i.e.
\[
- \pi(f, W f)_B +\pi \int_{B}(1-r^2)f^2dx  =  \pi \int_{B}(1-r^2)|\nabla f|^2dx + \pi \int_B f^2dx \, .
\]
We conclude that $-\pi(f, W f)_B$ is an entropy quantity, i.e. a measure of information, that we call the {\it prolate entropy} of $f$ w.r.t. $B$. 
In other words, $-\pi E_B W $ is an entropy operator. The {\it lucky accident} \cite{SP1}, that $W$ commutes with the truncated Fourier transform, finds a conceptual clarification in this fact. 

Based on the ordering of eigenvalues result in \cite{SP1}, we then have
\[
\text{lower prolate entropy} \longleftrightarrow \text{higher concentration}
\]
where the concentration is both on space and in Fourier modes as above. This is intuitive since information is the opposite of entropy. The above correspondence holds in the one-dimensional case, and we expect it to hold in general. 

In other words, in order to maximize simultaneously both quantities in  \eqref{SPe}, we have to minimize the prolate entropy. 

\section{Higher-dimensional Legendre operator}
The Legendre operator is the one-dimensional linear differential operator $\frac d{dx}(1-x^2)\frac d{dx}$. It is a Sturm-Liouville operator, probably best known because its eigenfunctions on $L^2(-1, 1)$ are the Legendre polynomials. 
In the following, we consider a natural higher-dimensional generalization of this operator. 

Let $\S$ be the Schwartz space of smooth, rapidly decreasing functions, $d\geq 1$. For the moment, we deal with complex-valued functions; the corresponding results for real-valued functions are obtained by restriction. 
We denote by $L_{\rm min} $ the $d$-dimensional {\it Legendre operator}, acting on $\S$, that we define by 
\ben\label{L}
L_{\rm min} =  \nabla (1-r^2)\nabla \, ;
\een
namely, $L_{\rm min} $ is the divergence of the vector field $(1-r^2)\nabla$, where $\nabla$ denotes the gradient and $r$ the radial coordinate in $\RR^d$. 
$L_{\rm min} $ can be written as
\ben\label{L2}
L_{\rm min}  = (1 -r^2)\nabla^2 - 2r\partial_r  \, ,
\een
indeed
\[
\nabla  (1-r^2)\nabla = \sum_k \partial_k\big((1- r^2)\partial_k   \big)
= \sum_k -2 x_k  \partial_k + \sum_k (1- r^2)\partial^2_k   
= -2r\partial_r  + (1 - r^2)\nabla^2 \, .
\]
We consider $L_{\rm min} $ as a linear operator on the Hilbert space $L^2(\RR^d)$, with domain $D(L_{\rm min} ) = \S$. The quadratic form associated with $L_{\rm min} $ is
\ben\label{Lher}
(f, L_{\rm min}  g) = -\int_{\RR^d}(1 - r^2) \nabla \bar f\!\cdot\!\nabla g\, dx\, ,\quad f,g\in \S\, ,
\een
because, by integration by parts, we have 
\begin{multline}\label{quad}
(f, \nabla\! \cdot \! (1-r^2)\nabla g) = \sum_k (f, \partial_k [(1-r^2)\partial_k g]) \\ = -\sum_k (\partial_k f,  [(1-r^2)\partial_k g])
= -\int_{\RR^d}(1 - r^2) \nabla g\!\cdot\!\nabla \bar f\, dx \, .
\end{multline}
\blem\label{Lemher}
 $L_{\rm min} $ is a Hermitian operator. 
\elem
\proof
Equation \eqref{Lher} shows that
\[
(f, L_{\rm min}  g) = ( L_{\rm min}  f,   g) \, ,
\]
for all $f,g\in\S$, therefore $L$ is Hermitian. 
\endproof
Thus  $L_{\rm min}\subset L_{\rm max}$, where $L_{\rm max} \equiv L^*$ denotes the adjoint of $L_{\rm min}$.    
\blem\label{L*}
$D(L_{\rm max})$ is the set of all $f\in L^2(\RR^d)$ such that $\nabla (1-r^2)\nabla f \in L^2(\RR^d)$ in the distributional sense, and $L_{\rm max} f =  \nabla  (1-r^2)\nabla f $ on $D(L_{\rm max})$. 
\elem
\proof
Let $f\in L^2(\RR^d)$, in particular $f\in S'(\RR^d)$ is a tempered distribution. With $g\in\S$, we have
\[
(f, \nabla  (1-r^2)\nabla g)  = \langle \nabla (1-r^2)\nabla f , g\rangle \, ,
\]
where the latter means the value of the distribution $\nabla (1-r^2)\nabla f$ on the test function $g$. 
Now, $f\in D(L_{\rm max})$ iff the linear functional
$g\in\S \mapsto (f, \nabla  (1-r^2)\nabla g)$ is continuous on $L^2(\RR^d)$, therefore iff $\nabla (1-r^2)\nabla f\in L^2(\RR^d)$ by Riesz lemma.
\endproof
Let $B$ be the unit open ball in $\RR^d$ and $E_B$ the orthogonal projection of $L^2(\RR^d)$ onto $L^2(B)$, that is $E_B$ is the multiplication operator by the characteristic function $\chi_B$ of $B$. Note that
\[
(f, L f) \leq 0\, ,\quad f\in \S, \ {\rm supp}(f)\subset \bar B\, ,
\]
as follows from \eqref{Lher}. 
\blem\label{sL}
Let $f,g$ be smooth functions on $\RR^d$. We have
\ben\label{LL}
\int_B  f \nabla (1-r^2)\nabla g  =  -\int_B (1 -r^2) \nabla  f\!\cdot \!\nabla g \, .
\een 
\elem
\proof
Taking into account that the vector field  $G = (1-r^2)\nabla g$ vanishes on $\partial B$, we have
\[
\int_B f \nabla \big((1-r^2) \nabla g\big ) 
  = \int_B f {\rm div}\, G  = - \int_B G\! \cdot\! \nabla f + \int_{\partial B} f G\!\cdot\! {\mathbf n}
  = - \int_B (1 -r^2) \nabla g\! \cdot\! \nabla f \, ,
\]
thus \eqref{LL} holds. 
\endproof
$L_{\rm min}$ does not commute with $E_B$, however, the following holds. 
\bprop\label{LE}
Let $f\in \S$. Then $\chi_B f \in D(L_{\rm max})$ and we have 
\ben\label{Lq}
L_{\rm max}\chi_B f = \chi_B L_{\rm min} f\, .
\een
Moreover, $L_{\rm max}$ is Hermitian on $\S + \chi_B \S$.  
\eprop
\proof
To prove the first part of the statement, namely eq. \eqref{Lq},
we must check that, for every $g\in\S$, we have $( f, \chi_B L_{\rm min} g) = (\chi_B  L_{\rm min} f, g)$, that is
\ben\label{her}
(\chi_B f, L_{\rm min} g) = ( L_{\rm min} f, \chi_B  g)\, .
\een
Taking into account that the vector field  $G = (1-r^2)\nabla g$ vanishes on $\partial B$, by eq. \eqref{LL} we have
\ben\label{LqB}
(\chi_B f, L_{\rm min} g) = \int_B\bar f L_{\rm min} g  = \int_B \bar f \nabla \big((1-r^2) \nabla g\big ) 
    = - \int_B (1 -r^2) \nabla \bar f\!\cdot \!\nabla g \, ,
\een
thus \eqref{her} holds because the last term in the above equality is symmetric in $f$ and $g$. 
\endproof
We shall denote by $L$ the closure of the restriction of $L_{\rm max}$ to $\S + \chi_B \S$. By Prop. \ref{Lq}, $L$ is Hermitian and commutes with $E_B$. 

Given $f\in L^2(\RR^d)$, we denote by $\hat f$ its Fourier transform 
\[
\hat f(p) = (2\pi)^{-d/2}\int_{\RR^d} e^{-ix\cdot p} f(x)dx\, ,
\] 
and by $\F$ the Fourier transform operator: $\F f =\hat f$. By Plancherel theorem, $\F$ is a unitary operator on $L^2(\RR^d)$.

In Fourier transform, $L_{\rm min}$ is given by the operator $\hat L_{\rm min} = \F L_{\rm min} \F^{-1}$; clearly $D(\hat L_{\rm min}) = \S$. 
We denote by
\ben\label{M}
M = (1- r^2)
\een
the multiplication operator by $(1- r^2)$ on $L^2(\RR^d)$. 
\blem
$\hat L_{\rm min} = - {r}^2 (1+\nabla^2) - 2r\partial_r $ on $\S$, where $r$  denotes the radial coordinate $|p|$ also in the dual space $\RR^d$. Therefore
\ben\label{LF1}
\hat L_{\rm min}    = L_{\rm min} - (\nabla^2 +1) + M\, .
\een
\elem
\proof
With $f\in\S$, we have
\[
-\big( (1 - r^2)\nabla^2 f\big){}^{\widehat{}}\,(p) = (1+\nabla_p^2)(|p|^2\hat f) 
= |p|^2\hat f + 2d \hat f + |p|^2\nabla^2_p\hat f + 4p\cdot\nabla_p\hat f
\] 
therefore, taking into account the equality $p\cdot\nabla_p = r\partial_r$, 
\[
\F\big( (1 - r^2)\nabla^2 \big)\F^{-1} =  - r^2(1 + \nabla^2) - 4r\partial_r - 2 d\, .
\]
On the other hand,
\[
\F( r\partial_r )\F^{-1} = -  r\partial_r  - d\ ,
\]
hence, accordingly with the expression \eqref{L2},
\[
\hat L_{\rm min}  = \F\big( (1 -r^2)\nabla^2 - 2r\partial_r \big)\F^{-1} = - r^2(1 + \nabla^2) - 2r\partial_r \, .
\]
Therefore
\ben\label{LF}
\hat L_{\rm min}   = (1 -r^2)\nabla^2 - 2 r\partial_r  
- (\nabla^2 +1) + (1- r^2)  = L_{\rm min} - (\nabla^2 +1) + M\, .
\een
\endproof
\section{Higher-dimensional prolate operator}
\label{Prol}
We now extend to the higher dimension some results in \cite[Sect. 1]{CMo}. 

Let $W_{\rm min}$ be the operator on $L^2(\RR^d)$ given by
\ben\label{WL0}
W_{\rm min} =  \nabla (1 -r^2)\nabla    - r^2   =  L_{\rm min} - r^2 
\een
with $D(W_{\rm min}) = \S$. $W_{\rm min}$ is a higher-dimensional generalisation of the {\it prolate operator}. 

By Prop. \ref{Lemher}, $W_{\rm min}$ is a Hermitian, being a Hermitian perturbation of $L_{\rm min}$ on $\S$; moreover, 
\[
- W_{\rm min} \geq - L_{\rm min} \geq 0
\]
on $D(W_{\rm min})\cap L^2(B)$, so $-W_{\rm min}$ is a positive operator on this domain. 

We explicitly note the equality
\ben\label{WL2}
-L_{\rm min}  =  -W_{\rm min} + M - 1   
\een
on $\S$ and that
\ben\label{LP}
- L_{\rm min}  \leq - W_{\rm min} \leq - L_{\rm min}+ 1 \quad\text{on $L^2(B)\cap D(L_{\rm min})$}\, ,
\een
because $0 \leq M \leq 1$ on $L^2(B)$. 
\bprop\label{WF}
$W_{\rm min}$
commutes with the Fourier transformation $\F$:
\[
\widehat{W}_{\rm min} = {W}_{\rm min}\, .
\]
Any linear combination of $L_{\rm min}$ and $M$ commuting with $\F$ is proportional to $W_{\rm min}$. 
\eprop
\proof
We have $\hat M = 1 +\nabla^2 $, therefore \eqref{LF1} gives
$\hat L_{\rm min}  = L_{\rm min} + M - \hat M$,
thus 
\[
L_{\rm min}+ M = \hat L_{\rm min} + \hat M
\]
on $\S$. By \eqref{WL2}, we then have 
\ben\label{WL}
W_{\rm min} = 
L_{\rm min} + M -1\, ,
\een
so
$W_{\rm min}  = \F W_{\rm min}  \F^{-1}$, 
as desired. 

Finally, if $a\in\RR$, we have
\[
\F(L_{\rm min} + a M)\F^{-1} = (L_{\rm min} + M - \hat M) + a\hat M = (L_{\rm min} + a M) +(1-a)(M -\hat M)\, ,
\]
thus $L_{\rm min} + aM$ commutes with $\F$ iff $(1-a)(M -\hat M) = 0$, that is iff $a=1$. 
\endproof
Let $\hat E_B = \F E_B \F^{-1}$ be the Fourier transform conjugate of the orthogonal projection $E_B:L^2(\RR^d)\to L^2(B)$, thus 
$(\hat E_B f)\hat{} = \chi_B \hat f$. In other words,
\[
\hat E_B f ={{(2\pi)}^{-\frac{d}2}}\tilde \chi_B * f \, ,
\]
where tilde denotes the Fourier anti-transform and $*$ the convolution product. We put $W_{\rm max} = W^*_{\rm min}$. We have
\[
D(W_{\rm max}) = \big\{f \in L^2(\RR^d) : \nabla (1 -r^2)\nabla f   - r^2 f \in L^2(\RR^d) \ \text{(distributional sense)} \big\}
\]
and 
\ben\label{Wtd}
W_{\rm max} f =  \nabla (1 -r^2)\nabla f   - r^2 f\, ,\quad f\in D(W_{\rm max})\, ,
\een
in the distributional sense.
Clearly, by Prop. \ref{WF}, also $W_{\rm max}$ commutes with $\F$
\ben\label{WFm}
W_{\rm max} = \F W_{\rm max} \F^{-1}\, .
\een
\bprop\label{WE}
Let $f\in\S$. Then $E_B f, \hat E_B f \in D(W_{\rm max} )$ and  
\ben\label{WEe}
W_{\rm max} E_B f =E_B  W_{\rm min}  f, \quad W_{\rm max} \hat E_B f = \hat E_B  W_{\rm min}  f \, .
\een
\eprop
\proof
Clearly $M$ commutes with $E_B$. Since $W_{\rm min}  = L_{\rm min} + M - 1$ \eqref{WL}, it follows from Prop.  \ref{LE} that
 $E_B f \in D(W_{\rm max})$ and $W_{\rm max}  E_B f =E_B  W_{\rm min}  f$, namely the first equation in \eqref{WEe} holds. 
 
 The second equation then follows from the first one by applying the Fourier transform because $W_{\rm min}$, $W_{\rm max}$ commute with $\F$, $\hat E_B = \F E_B \F^{-1}$, and $\F \S= \S$. 
\endproof
By the above proposition, we have
\[
\D \equiv \S + \chi_B\S + \widehat{\chi_B\S} \subset D(W_{\rm max} )
\]
and
\ben\label{Rm}
W_{\rm max} ( f + \chi_B g + \hat\chi_B * h ) = W_{\rm min}  f + \chi_B W_{\rm min} g + \hat\chi_B * W_{\rm min} h\, , \quad f,g,h\in \S\, ;
\een
recall that $\hat\chi_B$ is a smooth $L^2$-function vanishing at infinity, $\hat\chi_B(p) = \sqrt{\frac2{\pi}}\frac{\sin p}{p}$ if $d=1$. 
\blem\label{loc}
Let $f\in D(W_{\rm max})$ be a smooth function. Then, also the function $\chi_B f\in D(W_{\rm max})$, 
and $W_{\rm max}\chi_B f = \chi_BW_{\rm max} f$.
\elem
\proof
If $f\in\S$ the lemma follows as in Prop. \ref{WE}. Let now $f\in D(W_{\rm max})$ be a smooth function. Choose $f_0\in\S$ that is equal to $f$ on a neighborhood of $\bar B$. Then $\chi_B f = \chi_B f_0$, so $\chi_B f\in D(W_{\rm max})$. Moreover, 
\[
W_{\rm max} \chi_B f = W_{\rm max} \chi_B f_0 =  \chi_B W_{\rm min} f_0  =  \chi_B W_{\rm max} f\, ,
\]
where the last equality follows because $W_{\rm max}$ acts locally on $f$ by \eqref{Wtd}, 
so $W_{\rm max} f = W_{\rm min} f_0$ on a neighbourhood of $\bar B$.
\endproof
\blem\label{gdom}
For every $g\in\S$, $E_B  \hat E_B g$ belongs to $D(W_{\rm max})$ and we have
\ben\label{EE}
E_B W_{\rm max} \hat E_B g = W_{\rm max}  E_B  \hat E_B g\, .
\een
\elem
\proof
We may apply Lemma \ref{loc} with $f =  \hat E_B g$; indeed 
 $f = \hat\chi_B * g$ is a smooth function because $g$ is smooth, and $f$ in the domain of $W_{\rm max}$ by Prop. \ref{WE}.  
\endproof
Recall that a closed linear operator $Z$ on a Hilbert space $\H$ commutes with the orthogonal projection $F$ on $\H$ if
\ben\label{comm}
ZF \supset FZ\, ;
\een
this means 
\[
u\in D(Z) \implies Fu \in D(Z) \ \& \ ZFu = FZu\, .
\]
If $\D\subset D(Z)$ is a core for $Z$, then it suffices to verify the above condition for all $u\in\D$. 

Denote by $\F_B =  E_B\F E_B$ the {\it truncated Fourier transform}. Note that
\[
\F_B^*\F_B = E_B\F^* E_B\F E_B = E_B \hat E_B E_B
\]
is the {\it angle operator}. 
\bprop
\label{PropW}
The restriction of $W_{\rm max}$ to $\D$ is Hermitian. 
Its closure 
\[
W = \ov{W_{\rm max}|_\D}
\]
is Hermitian and
commutes with $\F$ and $E_B$, thus with $\hat E_B$ and $\F_B$ too. 
\eprop
\proof
$W_{\rm max}|_\D$ commutes with $\F$ because $W_{\rm max}$ commutes with $\F$, and $\D$ is globally $\F$-invariant. 

We now show that $W_{\rm max}|_\D$ is Hermitian. First, note that
$W_{\rm max}$ is Hermitian on $\S + \chi_B\S$ by the eq. \eqref{WL}, because $L_{\rm max}$ is Hermitian on $\S + \chi_B\S$ by Prop. \ref{LE}, and $M$ is Hermitian too on this domain. It then follows that $W_{\rm max}$ is Hermitian on $\S + \widehat{\chi_B\S}$ too due to \eqref{WFm}. 

So we have to show that $W_{\rm max}$ is symmetric on mixed terms in \eqref{Rm}. By \eqref{LL}, we are indeed left to check that  
$(\hat E_B g,  W_{\rm max} E_B h) = (W_{\rm max}\hat E_B g,  E_B h)$, for all $g,h\in\S$. 

Now,  by  Prop. \ref{WE} and Lemma \ref{gdom}, we have
\begin{multline*}
( \hat  E_B g,  W_{\rm max} E_B h) = (\hat  E_B g, E_B W_{\rm min} h) = ( E_B \hat  E_B g, W_{\rm min} h) 
= (W_{\rm max}  E_B  \hat E_B g,  h) \\
= (E_B W_{\rm max} \hat E_B g,  h)= (E_B  \hat E_B W_{\rm min} g,  h) = (  \hat E_B W_{\rm min} g, E_B h) 
= ( W_{\rm max} \hat E_B g,  E_B h) \, .
\end{multline*}
So $W_{\rm max}|_\D$ is Hermitian, hence its closure $W$ is Hermitian too. 

It remains to show that $W $ commutes with $E_B$. 
We need to check that 
$E_B W  \subset  W  E_B$. With $f,g,h\in\S$, we then have to verify that $E_B( f + E_B g + \hat E_B h )$ belongs to the domain of $W $ and
\ben\label{WEB}
 W  E_B( f + E_B g + \hat E_B h ) = E_B  W  ( f + E_B g + \hat E_B h )  \, .
\een
By linearity, we can check the above condition for each of the three terms individually. Concerning the first term, that is the case $g=h=0$, we have $f\in \D \subset D(W )$ and by Prop. \ref{WE}
\[
W  E_B f = W_{\rm max} E_B f =  E_B W_{\rm min} f = E_B W  f\, .
\]
Consider now the last term. With $k \equiv \hat E_B h$, we have to show that $E_B k$
belongs to $D(W )$ and $W  E_B k = E_B W  k$. As $k$ is smooth, we can pick $k_0\in\S$ that agrees with $k$ in a neighborhood of $\bar B$; so $E_B k = E_B k_0\in D(W )$, and  by Lemma \ref{loc}  we have
\[
W  E_B k = W_{\rm max} E_B k =  W_{\rm max} E_B k_0 = E_B W_{\rm max} k_0 =  E_B W_{\rm max} k = E_B W  k\, ,
\]
where the equality $E_B W_{\rm max} k_0 =  E_B W_{\rm max} k$ follows  as in the proof of Lemma \ref{loc} because $W_{\rm max}$ acts locally. 

Concerning the remaining second-term case, take then $g\in\S$; clearly $E_B E_B g = E_B g \in\D \subset D(W )$ and we are left to show that $W  E_B g = E_B W  E_B g$.  

Now, $W_{\rm max} E_B g$ is supported in $\bar B$ as distribution because $ W_{\rm max}$ is local; on the other hand,  
$W_{\rm max} E_B g$ is an $L^2$-function, therefore $W_{\rm max} E_B g =  E_B W_{\rm max} E_B g$. 
We conclude that
\[
W  E_B g =  W_{\rm max} E_B g =  E_B W_{\rm max} E_B g =  E_B W  E_B g \, ,
\] 
and the proof is complete. 
\endproof
$W$ is the minimal closed extension of $W_{\rm min}$ that commutes both with $E_B$ and $\hat E_B$. Indeed, if $\tilde W$ is an extension of $W_{\rm min}$ with this property, then $D(\tilde W)$ must contain $E_B\D$ and $\tilde WE_B f = E_B W_{\rm min} f$, $f\in \S$. Similarly with $\hat E_B$ in place of $E_B$. So $\tilde W \supset W$. 

Note that the angle operator $E_B \hat E_B E_B$ is of trace class, indeed $E_B \hat E_B |_{L^2(B)}$ is the positive Hilbert-Schmidt $T_B$ on $L^2(B)$ operator with kernel $k_B(x-y)$ where
\[
k_B({ z}) = \frac1{(2\pi)^{d/2}}\int_B e^{-i x\cdot { z}}dx\, \chi_B({z})
\]
($k_B = \hat \chi_B$ on $B$, zero out of $B$). 
The eigenvalues of $T_B$  are strictly positive, $\l_1 > \l_2 >\cdots \l_k>\cdots >0$, with finite multiplicity. 
The equality
\[
||\F f||^2 = (f, \F_B^*\F_B f) = (f_B, T_B f_B)_B\, ,
\]
$f_B = f|_B$, shows that the normalized $k$-th eigenfunctions of $T_B$ are concentrated at level $\l_k$ in an appropriate sense. 
Note that, on the even function subspace, $\F$ is a unitary involution, thus $\F_B$ is selfadjoint; so $\F_B$ and $\F_B^* \F_B = E_B\hat E_B E_B$ share the same eigenfunctions.  

We now show that $- W $ is positive on $B$, namely $-E_B W$ is positive. 
\bprop\label{WsPp}
For every $u\in D(W )$, we have 
\ben\label{WsP}
- (u, W  u)_B = - \int_B\bar u  W  u\, dx \geq 0\, .
\een
\eprop
\proof
As $\D$ is a core for $W $, it suffices to check \eqref{WsP} with $u = f + E_B g + \hat E_B h$, with $f,g,h\in\S$. 

Now, $\chi_B u$ is a smooth function on $\bar B$; choose $u_0\in\S$ that agrees with $u$ on $\bar B$. 
By eq. \eqref{WEB}, we have
\[
\int_B \bar u W  u dx = \int_B \bar u W_{\rm max} u dx 
= \int_B  \bar u_0 W_{\rm min}  u_0 dx =- \int_B \bar u_0 L u_0 dx  -  \int_B |x|^2  |u_0|^2 dx \leq 0
\]
by \eqref{LqB}, because $W_{\rm max}$ is local. 
\endproof
As seen, both $W$ and $L$ commute with $E_B$, and we consider now their restrictions to $L^2(B)$, which we denote by $W_B$ and $L_B$. 

Let $C^\infty(\bar B)$ be the space of smooth function on $\bar B$, up to the boundary; we may regard $C^\infty(\bar B)$ as a subspace of $L^2(B)\subset L^2(\RR^d)$. As is known, $\chi_B \S = C^\infty(\bar B)$. 
We now show that  $W_B$ and $L_B$ are essentially selfadjoint on  $C^\infty(\bar B)$. We will also denote by $C_0^\infty(B)$ the space of smooth functions on $\bar B$ with compact support contained in $B$.

\bcor\label{LWB}
Both $W_B$ and $L_B$ are selfadjoint, positive operators on $L^2(B)$. $C^\infty(\bar B)$ is a core for both $W_B$ and $L_B$. 
\ecor
\proof
As $W_B$ is Hermitian and commutes with the positive Hilbert-Schmidt operator $T_B$, it follows that $W_B$ is selfadjoint.  

Since $W$ commutes with $E_B$, $\D$ is a core of $W$ and $E_B \D \subset \D$, it follows that $E_B\D$ is a core for $W_B$. 
On the other hand, $E_B\D = \chi_B\S$ because functions in $\S + \hat E_B \S$ are smooth; so $\chi_B \D = C^\infty(\bar B)$. Therefore $ C^\infty(\bar B)$ is a core for $W_B$. 
$W_B$ is then positive by Prop. \ref{WsPp}. 

Since $L_B$ is a bounded perturbation of $W_B$ on $L^2(B)$, also $L_B$ is selfadjoint with core $C^\infty(\bar B)$. $L_B$ is then positive by  Lemma \ref{sL}. 
\endproof
In the one-dimensional case, the essentially selfadjointness of $L_B$ on $C^\infty[-1,1]$ (thus of its bounded perturbation $W_B$) follows by the well-known fact that the Legendre polynomials form a complete orthogonal family of  $L_B$-eigenfunctions.  Note that $L_B$ is not essentially selfadjoint on $C_0^\infty(-1,1)$, see \cite{K}. 
\bprop\label{Fri}
$C_0^\infty(B)$ is a form core for $L_B$, thus for $W_B$. Moreover, $-L_B$ and $-W_B$ are the Friedrichs extensions 
of $-L_B|_{C_0^\infty(B)}$ and $-W_B|_{C_0^\infty(B)}$. 
\eprop
\proof
We consider $L_B$ only because $W_B$ is a bounded perturbation of it. Since $L_B$ is essentially selfadjoint on $C^\infty(\bar B)$, it is enough to show that the form closure of the quadratic form $q$ of $-L_B|_{C_0^\infty(B)}$ contains $C^\infty(\bar B)$. 

Now, $q$ is given by \eqref{LL} on $C_0^\infty(B)$. By \cite[Prop. 10.1]{Sch}, it suffices to show that, given $u\in C^\infty(\bar B)$, there exists a sequence of functions $u_n\in C_0^\infty(B)$ such that $u_n\to u$ and 
\[
q(u_n, u_n) = \int_B(1 -r^2)|\nabla u_n|^2dx
\] 
is bounded, $n\in \mathbb N$. First suppose $u = \chi_B$. Let $h_n\in C_0^\infty(-1,1)$ be even such that $h_n = 1$ on $(0,1 - \frac1{n})$
and $|h'_n|$ bounded by $2 n$ and set $u_n(x) = h_n(r)$. Then $u_n \to \chi_B$ and the sequence
\[
q(u_n, u_n) =  \int_B(1 -r^2)|\nabla u_n|^2dx 
\leq \int_{1 - 1/n\leq r\leq 1}(1 -r^2)(2n)^2dx 
\leq {\rm const.} \frac1{n^2}
(2n)^2
\]
is bounded. The case of a general $u\in C^\infty(\bar B)$ follows on the same lines by replacing $u_n$ by $u_n u$. 

So, $C^\infty(\bar B)$ is in the domain of the square root $\sqrt{-L_F}$ of the Friedrichs extension $-L_F$ of $-L_B|_{C_0^\infty(B)}$. On the other hand,
$C^\infty(\bar B)$ is a core for $-L_B$, thus for $\sqrt{-L_B}$. We conclude that $L_B = L_F$. 
\endproof
See e.g. \cite{Sch} for the Friedrichs extension. 

\section{Modular theory and entropy of a vector}
In this section, we recall the basic structure concerning the modular theory of a standard subspace $H$, the entropy of a vector relative to $H$, and their applications to the entropy density of a wave packet. 

\subsection{Entropy operators}
\label{entroper}
Let $\H$ be a complex Hilbert space and $H\subset \H$ a standard subspace, i.e. $H$ is a real linear, closed subspace of $\H$ such that $H\cap iH = \{0\}$ and ${\ov{H + iH}} = \H$, with $H'$ the symplectic complement of $H$, 
\[
H' = \{\Phi'\in \H : \Im(\Phi,\Phi')=0,\, \Phi\in H\}\, . 
\]
The Tomita operator \[
S_H: \Phi_1 + i\Phi_2\in H + iH \mapsto \Phi_1 - i\Phi_2\in H+ iH\, ,\quad \Phi_1,\Phi_2\in H\, ,
\]
is anti-linear, closed, densely defined, and involutive on $\H$. Let $S_H = J_H \Delta_H^{1/2}$ be the polar decomposition of $S_H$.  $\Delta_H$ is called the {\it modular operator} associated with $H$; it is a canonical positive, non-singular selfadjoint operator on $\H$ that satisfies 
\[
\Delta^{is}_H H = H\, , \quad s\in \mathbb R\, .
\]
The one-parameter unitary group $s\mapsto \Delta^{is}_H$ on $\H$ is called the {\it modular unitary group} of $H$, whose generator $\log\Delta_H$ is the {\it modular Hamiltonian}. $J_H$ is an anti-unitary involution on $\H$ and $J_H H = H'$, named the {\it modular conjugation} of $H$. 

For simplicity, let us assume that $H$ is factorial, namely $H\cap H' =\{0\}$, see e.g. \cite[Sect. 2.1]{CLR}  for the general case of a closed, real linear subspace.   

The {\it entropy of a vector} $\Phi\in\H$ with respect to a standard subspace $H\subset\H$ is defined by
\ben
S_\Phi = S^H_\Phi = \Im(\Phi, P_H A_H\, \Phi) = (\Phi, P^*_H \log\Delta_H\, \Phi) 
\een
(in a quadratic form sense), where $P_H$ is the {\it cutting projection} 
\[
P_H : H + H' \to H\, , \quad \Phi +\Phi' \mapsto \Phi
\] 
and $A_H =  -i\log\Delta_H$ \cite{L19, CLR}, the semigroup generator $\frac{d}{ds}\Delta^{-is}_H |_{s=0}$ of the modular unitary group. 

We have $P^*_H = -iP_H i$
and the formula in \cite{CLR}
\ben\label{fP}
P_H = (1 - \Delta_H)^{-1} + J_H\Delta_H^{1/2} (1 - \Delta_H)^{-1}  \, ;
\een
($P_H$ is the closure of the right hand side of  \eqref{fP}).

The {\it entropy operator} $\E_H$ is defined  by
\ben\label{eh}
\E_H = iP_H i \log\Delta_H
\een
(closure of the right-hand side). We have
\ben
S_\Phi = (\Phi, \E_H \Phi) \, ,\quad \Phi\in \H\, .
\een
Here, $S_\Phi$ is defined for any vector $\Phi\in\H$ as follows. $S_\Phi = q(\Phi,\Phi)$ with $q$ the closure of the real quadratic form $\Re(\Phi,\E_H \Psi)$, $\Phi,\Psi\in D(\E_H)$. So $S_\Phi = +\infty$ if $\Phi$ is not in the domain of $q$. 
\bprop
The entropy operator $\E_H$ is real linear, positive, and selfadjoint w.r.t. to the real part of the scalar product.
\eprop
\proof
$\E_H$ is clearly real linear, and positive because the entropy of a vector is positive \cite[Prop. 2.5 (c)]{CLR}. The selfadjointness of $\E_H$ follows by the formula \eqref{fP},  see \cite[Lemma 2.3]{LM}. 
\endproof
In our view, an {\it entropy operator} $\E$ is a real linear operator on a real or complex Hilbert space $\H$, such $\E$ is positive, selfadjoint and its expectation values $(f,\E f)$, $f\in\H$, correspond to entropy quantities (w.r.t. $B$). $\E$ may be unbounded, and $(f,\E f)$ is understood in the quadratic form sense, so it takes values in $[0,\infty]$. 
It is convenient to consider more entropy operators by performing operations, that preserve our demand, on the entropy operators.  

\smallskip\noindent
{\it Basic.} If $\E$ is a real linear operator on a real Hilbert space $H$ of the form \eqref{eh}, we say that $\E$ is an entropy operator. 

\smallskip\noindent
{\it Restriction and direct sum.} If $\E = \E_+\oplus \E_-$ on a {\it real} Hilbert space direct sum $H = H_+\oplus H_-$, then $\E$ is an entropy operator on $H$, iff both $\E_\pm$ are entropy operators.

\smallskip\noindent
{\it Change of metric.} Suppose that ${\cal S} \subset H$ is a core for the entropy $\E$ on $H$ and $(\cdot,\cdot)'$ is a scalar product on ${\cal S}$; denote by $H'$ the corresponding real Hilbert space completion and by $\jmath: {\cal S}\subset H' \to H$ the identification map. If $\jmath^* \E \jmath$ is densely defined, its Friedrichs extension $\E'$ is an entropy operator on $H'$. Note that
\[
(f, \E' f)' = (f, \E f)\, ,\quad f\in {\cal S}\, .
\]
\smallskip\noindent
{\it Sum, difference.} If $\E_1,\E_2$ are entropy operators and $\E = \E_1 \pm \E_2$ is densely defined and positive, the Friedrichs extension  $\E$ is an entropy operator. 

\smallskip\noindent
{\it Born entropy.} $\pi E_B$, with $E_B$ the orthogonal projection onto $L^2(B)$,  is an entropy operator on $L^2(\RR^d)$. 

\smallskip\noindent
In order to justify the last item, note that $(f, E_B f) = ||f||^2_B$. In Quantum Mechanics, with the normalization $|| f ||^2 = 1$, $|| f ||^2_B$ is the particle probability to be localized in $\Omega$, accordingly to Born's interpretation. Moreover, in Communication Theory, $||f||^2_B$ represents the part of energy of $f$ contained in $B$ \cite{SP1}. We thus define 
\ben\label{Born}
\pi(f, E_B f) = \pi ||f||^2_B =  \pi \int_B f^2dx = \text{{\it Born entropy} of $f$ in $B$}\, ,
\een
$f\in L^2(\RR^d)$ real. 
The $\pi$ normalization is chosen by compatibility reasons (Sect. \ref{proe});

\subsection{Abstract field/momentum entropy}
We consider two real linear spaces ${\cal S}_+$ and ${\cal S}_-$ and a duality $f , g\in {\cal S}_+  \times {\cal S}_-\mapsto \langle f, g \rangle  \in \RR$.   A real linear, invertible operator 
\[
\mu : {\cal S}_+ \to {\cal S}_-
\]
is also given; we assume that $\mu$ is symmetric and positive with respect to the duality, i.e.
\ben\label{md}
\langle f_1, \mu f_2\rangle =  \langle f_2, \mu f_1 \rangle \, , \quad f_1 , f_2\in {\cal S}_+\, ,
\een
\[
\langle f, \mu f\rangle \geq 0\, , \quad f\in {\cal S}_+\, ,
\]
with $\langle f, \mu f\rangle = 0$ only if $f=0$. 

So ${\cal S}_\pm$ are real pre-Hilbert spaces with scalar products 
\[
(f_1, f_2)_+ = \langle f_1, \mu f_2\rangle \, , \ \  (g_1, g_2)_- = \langle \mu^{-1} g_2,  g_1 \rangle\, , \quad f_1, f_2\in {\cal S}_+,\  g_1, g_2\in {\cal S}_-\, ,
\]
and $\mu$ is a unitary operator. 

Let $H_\pm$ be the real Hilbert space completion of ${\cal S}_\pm$. Then $\mu$ extends to a unitary operator $H_+\to H_-$, still denoted by $\mu$. Moreover, the duality between ${\cal S}_+$ and ${\cal S}_-$ extends to a duality between $H_+$ and $H_-$
\[
\langle f, g\rangle  = (f, \mu^{-1}g)_+ = (\mu f, g)_-\, , \quad f\in H_+, g\in H_-\, .
\]
Set $\H = H_+ \oplus H_-$. The bilinear form $\b$ 
on $\H$
\ben\label{beta}
\b(\Phi, \Psi) =    \langle g_1, f_2\rangle -\langle f_1, g_2 \rangle
\een
$\Phi \equiv f_1\oplus g_1$, 
$\Psi \equiv f_2\oplus g_2$,
is symplectic and non-degenerate (the coefficient $\frac12$ is to conform with the next section case). 
This will be the imaginary part of the complex scalar product of $\H$:
$
\Im( \Phi ,  \Psi ) = \b( \Phi ,  \Psi )$.

Now, the operator
\ben\label{imum}
\imath = \left[\begin{matrix}
0 & \mu^{-1} \\ -\mu  & 0  
\end{matrix}\right]\, ,
\een
namely $\imath :  f\oplus g \mapsto \mu^{-1} g\oplus -\mu  f$, is a unitary on $\H =  H_+\oplus H_+$.

By \eqref{md}, $\imath$ preserves $\b$, that is $\b(\imath\Phi , \imath\Psi) = \b(\Phi ,\Psi)$. 
As $\imath^2 = -1$,  the unitary $\imath$ defines a complex structure (multiplication by the imaginary unit)  on
$\H$ that becomes a complex Hilbert space with a scalar product
\[
(\Phi,\Psi) = \b(\Phi, \imath \Psi) + i\b(\Phi, \Psi) 
\]
($i= \sqrt{-1}$).
That is
\[
( \Phi,  \Psi) = \big[\langle f_1,\mu f_2\rangle +\langle\mu^{-1}g_2 ,g_1 \rangle\big] + i \big[\langle f_2, g_1\rangle -\langle f_1 , g_2\rangle\big]\, ,
\]
$\Phi \equiv f_1\oplus g_1$, $\Psi \equiv f_2\oplus g_2$ as above.  

Suppose now $K_\pm \subset H_\pm$ are closed, real linear subspaces. The symplectic complement $K'$ of $K \equiv K_+ \oplus K_-$ is
\[
K' = \big\{ f\oplus g \in H : \b( f\oplus g , h\oplus k) = 0,\,  h\oplus k\in K\big\} = K_-^o\oplus K_+^o \, ,
\]
where $K_\pm^o$ denote the annihilatotors of $K_\pm$ in $H_\mp$ under the duality $\langle\cdot ,\cdot\rangle$. 

Let us consider the case $K$ is standard and factorial. Then the cutting projection
\[
P_K = K + K' \to K
\]
is diagonal
\[
P_K = \left[\begin{matrix}
P_+ & 0 \\ 0  & P_-  
\end{matrix}\right]\, ,
\]
with $P_\pm$ the projection $P_\pm: K_\pm + K^o_\mp \to K_\pm$. 
\bprop
The modular Hamiltonian $\log\Delta_K$  and conjugation $J_K$ are diagonal; so
$A_K = -\imath\log \Delta_K$ is off-diagonal, that is 
\ben\label{AKd}
A_K  = \pi \left[
\begin{matrix}
0 & {\bf M} \\ {\bf L} & 0
\end{matrix}\right] \, ,
\een
with ${\bf M}$ and ${\bf L}$ operators $H_\pm\to H_\mp$. 

The entropy of $\Phi \equiv f\oplus g\in\H$ with respect to $K$ is given by
\ben\label{Sfa}
S_\Phi = -\pi  \langle f,  P_-{\bf L}f\rangle  + \pi \ \langle g,  P_+{\bf M} g\rangle \, .
\een
In particular, if $\Phi\in K$,
\[
S_\Phi = -\pi \langle f, {\bf L}f\rangle  + \pi \langle g,  {\bf M} g\rangle \, .
\]
\eprop
\proof
As $K = K_+ \oplus K_-$ and $\imath K = \mu^{-1} K_-  \oplus \mu K_+$ are direct sum subspaces,
the Tomita operator $S_K$ is clearly diagonal, and so is its adjoint $S^*_K$. The modular operator $\Delta_K = S_K^* S_K$ is thus diagonal. Since the logarithm function is real on $(0,\infty)$, by functional calculus the modular Hamiltonian $\log\Delta_K$ is diagonal too. Also $J_K$ is diagonal due to formula \eqref{fP}. 

So $A_K$ is off-diagonal because $\imath$ is off-diagonal and we may write $A_K$ as in \eqref{AKd}. We have
\ben
P_K A_K  = \pi \left[
\begin{matrix}
0 & P_ + {\bf M} \\ P_- {\bf L} & 0
\end{matrix}\right] \, ,
\een
thus the entropy of $\Phi$ is given by
\[
S_\Phi = \b(f\oplus g, P_K A_K\, f\oplus g) =  \pi \b(f\oplus g, P_+ {\bf M} g\oplus P_+{\bf L} f)
=  - \pi  \langle f,  P_-{\bf L}f\rangle  +  \pi \langle g,  P_+{\bf M} g\rangle  \, .
\]
\endproof
The fact that $\log\Delta_K$ is diagonal was shown in \cite{BCM}, based on the the formula $P_K- \imath P_K \imath = 2(1 - \Delta_K)^{-1}$, which follows from \eqref{fP}.  

The entropy operator is given by
\ben
\E_K  =  \pi \left[
\begin{matrix}
-\mu^{-1}P_- {\bf L} & 0 \\ 0 & \mu P_ + {\bf M}
\end{matrix}\right] \, .
\een
Note that, since $A_K$ is skew-selfadjoint and complex linear on $\H$, we have the relations
\ben\label{M*}
{\bf M}^* = -{\bf L} = \mu {\bf M}\mu\, .
\een
Clearly, 
\[
- \pi \langle f, P_-{\bf L}f\rangle = S_{f\oplus 0}\, ,\quad  \pi  \langle g,  P_+{\bf M} g\rangle = S_{0\oplus g}\, .
\]
We then define: 
\begin{align*}
-\ \pi  \langle f, P_-{\bf L}f\rangle\quad  &\text{{\it field entropy} of}\ f\in {\cal S}_+\ w.r.t. \ K_+\, ,
\\
 \pi  \langle g,  P_+{\bf M} g\rangle\quad  &\text{{\it momentum entropy} of}\ g\in {\cal S}_- \, w.r.t. \ K_-\, .
\end{align*}
(quadratic form sense). 
Note that only the duality, not the Hilbert space structure, enters directly into the definitions of the above entropies. 
\subsection{Local entropy of a wave packet}
The above structure concretely arises in the wave space context; namely, in the free, massless, one-particle space in Quantum Field Theory. 

Denote by $\Sr$ the real Schwarz space. 
As is known, if $f, g\in\Sr$, there is a unique smooth real function $\Phi(t,\bx)$ on $\mathbb R^{1+d}$ which is a solution of the wave equation 
\[
\square \Phi \equiv\partial_t^2 \Phi - \nabla^2_{x}\Phi  = 0
\]
(a wave packet or, briefly, a wave) with Cauchy data $\Phi |_{t =0} = f$, $\partial_t\Phi |_{t =0} = g$. 
We set $\Phi=w(f,g)$ and denote
by $\T$ the real linear space of these $\Phi$'s; we will often use the identification
\ben\label{st}
\Sr\oplus\Sr \longleftrightarrow \T\, , \qquad \ f\oplus g \longleftrightarrow w(f, g)\, .
\een  
We thus deal directly with $\Sr\oplus\Sr$ and consider the symplectic form on it
\ben\label{beta2}
\b(f_1\oplus g_1, f_2\oplus g_2) =
 ( g_1, f_2) - ( f_1, g_2 ) \, ,
\een
where the scalar product in \eqref{beta2} is the one in $L^2(\RR^d)$. 

We set ${\cal S}_+  = \Sr$, ${\cal S}_-  = \mu\Sr$, with $\mu$ the given by
\ben\label{mum}
\widehat{\mu f}({ p}) = | p|\hat f({ p})\, .
\een
The duality between ${\cal S}_+$ and ${\cal S}_-$ is given by the $L^2$ scalar product. 
Let $H_{\pm}$ be the real Hilbert space of tempered distributions $f \in \Sr'$ such that $\hat f$ is a Borel function with
\ben\label{H12}
 \|f \|^2_{\pm} = \int_{\RR^d} {|p|^{\pm1}} | \hat f(p)|^2 dp< +\infty\,  .
\een
${\cal S}_\pm$ is dense in $H_{\pm}$, yet $\S\subset H_-$ only if $d>1$. 
The complex Hilbert space is
$\H$
is the real Hilbert space $H  = H_+\oplus H_-$ equipped with complex structure given  $\imath$ \eqref{imum}. 

With
\[
H_{\pm}(B) = \big\{f_\pm \in {\cal S}_\pm : {\rm supp}(f_\pm) \subset B\big\}^-\, ,
\]
the standard subspace $K\equiv H(B)\subset\H$ is 
\[
H(B) = H_{+}(B)\oplus H_{-}(B)\, .
\]
Set $\Delta_B = \Delta_{H(B)}$ for the modular operator associated with $H(B)$, and $A_B = -\imath \log \Delta_B$. The action of $\Delta_B^{is}$, $s\in\RR$, on $\T$ is geometric \cite{HL}, so $A_B$ is computable. 
\bthm\label{thmA} {\rm \cite{LM}}.
 On $\Sr\times\Sr$, $d>1$, we have
\ben\label{Kmat}
A_B  = \pi \left[
\begin{matrix}
0 & (1 - r^2) \\ (1 - r^2)\nabla^2 - 2r\partial_r  - 2D & 0
\end{matrix}\right] \, ,
\een
namely 
\ben\label{AB}
A_B  = \pi \left[
\begin{matrix}
0 & M \\ L_D& 0
\end{matrix}\right]
\een
with $L_D = L - 2D$; 
here, $L: H_+ \to H_-$,
$M:H_-\to H_+$ are the closure of the operators \eqref{Le}, \eqref{M} on $\S$, and $D = (d-1)/2$ (the scaling dimension).  

Case $d=1$: the above formula still holds on $S_{\rm r}(\RR)\times{\dot S}_{\rm r}(\RR)$, with ${\dot S}_{\rm r}(\RR)$ the subspace of $S_{\rm r}(\RR)$ consisting of functions with zero mean \cite{L22}. 
\ethm
\noindent
In the following, we assume $d >1$. The case $d=1$ is similar, it is sufficient to replace $S_{\rm r}(\RR)\times S_{\rm r}(\RR)$ by 
$S_{\rm r}(\RR)\times{\dot S}_{\rm r}(\RR)$ as above. 
\bcor{\rm \cite{LM}.}
Let $\Phi = w(f,g)$ be a wave packet with Cauchy data $f,g\in \Sr$. The entropy of $\Phi$ in $B$ (i.e. with respect to $H(B)$) is given by
\[
S_\Phi = - \pi (f,  Lf)_B +  \pi  (g,  Mg)_B + 2\pi D ||f||^2_B 
\]
($L^2$-scalar product)\footnote{
The symplectic form in \cite{CLR,LM} is defined as $1/2$ the one given by \eqref{beta2}. The entropy values in this paper are  twice the ones  there.}
. 
\ecor
\proof
The corollary follows by \eqref{Sfa} because the cutting projection $P_{H(B)}$ is given by the multiplication by the characteristic function $\chi_B$ on both components of $\Sr\times\Sr$, and the duality $\langle\cdot , \cdot\rangle$ by the $L^2$ scalar product. So
\begin{multline*}
S_\Phi = - \pi ( f,  \chi_B {L_D}f)  + \pi (g,  \chi_B {M} g) \\
=  \pi\int_{B}(1-r^2)|\nabla f |^2 d\bx 
+ \pi D \int_B f^2dx +  \pi\int_{B}(1-r^2)g^2 dx\\
 =  - \pi (f, Lf)_B + \pi (g, Mg)_B + \pi D ||f||^2_B \ ,
\end{multline*}
by the equality \eqref{LqB}. 
\endproof
More generally, if $f,g\in L^2(\RR^d)$, we set
\ben\label{Sgen}
S_{f\oplus g} =  - \pi (f, L_B E_B f) +  \pi (g, Mg)_B + 2\pi D ||f||^2_B \ ,
\een
in the quadratic form sense. 
As a consequence, we have a lower bound for the entropy. 
\bcor
The entropy of $\Phi = f\oplus g$ in $B$, $f,g\in L^2(B)$, is lower bounded by
\ben\label{Slb}
S_\Phi \geq 2\pi D ||f||^2_B\, .
\een
The inequality \eqref{Slb} is an equality if $f = \chi_B$, $g =0$; in this case
\[
S_{f\oplus g}  = 2\pi {\rm Vol}(B) D\, .
\]
\ecor
\proof
The inequality \eqref{Slb} is immediate as both terms $- \pi (f, Lf)_B $ and $\pi (g, Mg)_B$ are non-negative. 

Since $\chi_B$ belongs to the domain of $L_B$  and $L_B\chi_B = 0$, the inequality  is an equality if $f = \chi_B$, $g =0$ by \eqref{Sgen}. 
\endproof
Note that, since $-L_{\rm min} =  -W_{\rm min} + M - 1$,  we may rewrite $S_\Phi$ as follows,
$\Phi = w(f,g)$: we have
\ben\label{SP}
S_\Phi = \pi \Big( -(f,  W  f)_B +  (f,  Mf)_B + (g,  Mg)_B + \frac{d-2}2 ||f||^2_B\Big)\, .
\een
We are going to see in the next section that each individual term on the right-hand side of the above equality has an entropy interpretation. 

We end this section by writing up the formula
$|\nabla| ( 1 - r^2 ) |\nabla| = - L + 2D$
on $\S$, which follows from \eqref{M*}, where $ |\nabla| = \sqrt{-\nabla^2}$, the square root of minus Laplacian on $L^2(\RR^d)$.

\section{Prolate entropy} 
\label{proe}
By Thm. \ref{thmA}, the modular Hamiltonian $\log \Delta_B = \imath A_B$ is the closure of the linear operator on $H=H_+\oplus H_-$ given by
\[
\log \Delta_B = \pi \left[
\begin{matrix}
-\mu L_D & 0 \\ 0 &\mu^{-1} M
\end{matrix}\right]
\]
with core domain $\S\oplus \S$. 

The cutting projection w.r.t. $H(B)$  is 
$P_B = \left[\begin{matrix}
\chi_B & 0 \\ 0 & \chi_B
\end{matrix}\right]$,
therefore the entropy operator $\E_B = \imath P_B \imath \log\Delta_B$ on $H$ is given by
\ben\label{e1}
\E_B = \left[
\begin{matrix}
- \pi \chi_B L_D & 0 \\ 0 &  \pi \chi_B M
\end{matrix}\right]\, .
\een
Let $\jmath_\pm: \S\subset L^2(\RR^d)\to  H_\pm$ be the identification map on $\S$. 
Then $\jmath^*_\pm = \mu^{\mp 1}$. 

The {\it entropy operator $\E'_B$ on $L^2(\RR^d)\oplus L^2(\RR^d)$} corresponding to $\E_B$ in the sense of Sect. \ref{entroper}  is therefore given by
\ben\label{e2}
\E'_B = \left[
\begin{matrix}
- \pi E_B L_D & 0 \\ 0 & \pi E_B M
\end{matrix}\right]\, .
\een
So each of the two components of $\E'_B$ is an entropy operator on $L^2(\RR^d)$; and so is $-E_B L_B = -E_B L_D - 2DE_B$, due to \eqref{Born}. 
 More precisely, $ME_B$ is essentially selfadjoint on $\S$;
the Friedrichs extensions of $-E_B L_B$  on $\S$ is equal to $L_B E_B$ by Prop. \ref{Fri}; 
so both $-\pi L_B E_B$ and $\pi M E_B$ are entropy operators on the Hilbert space $L^2(\RR^d)$ (see Sect. \ref{entroper}). 

With $f\in L^2(\RR^d)$ real, we set
\[
\pi (f,  M f)_B = \pi \int_B(1-r^2)f^2dx = \text{{\it parabolic entropy} of $f$ in $B$}\, .
\]
This is equal to the entropy $S_\Phi$ of the flat wave $\Phi = w(0,f)$. 

Similarly, we set
\[
-  \pi (f,  L f)_B = \pi \int_B(1-r^2)|\nabla f|^2dx = \text{{\it Legendre entropy} of $f$ in $B$}\, .
\]
This is equal to the entropy $S_\Psi$ of the stationary wave $\Psi = w(f,0)$. 

Now, $-LE_B  =  -WE_B + ME_B - E_B$, so $\pi WE_B$ is an entropy operator too; we thus define:
\[
- \pi (f,  W  f)_B  = \pi \int_B\big((1-r^2)|\nabla f|^2 + r^2\big)dx = \text{\it prolate entropy {\rm of}} \ f\ \text{in} \ B\, ,
\]
$f\in L^2(\RR^d)$ real. 

We summarize our discussion in the following theorem. 
\bthm\label{thme}
$- \pi W E_B$ is an entropy operator on $L^2(\RR^d)$. 
The sum of the prolate entropy and the parabolic entropy is equal to the sum of the Legendre entropy and the Born entropy, all with respect to $B$. Namely, the relation \eqref{ES} holds for every $f\in L^2(\RR^d)$, in the quadratic form sense. 

$WE_B$ commutes with the truncated Fourier transform $\F_B$. Let $V$ be a real linear combination of $LE_B , ME_B$ and $E_B$ commuting with $\F_B$; then $V = a WE_B + b E_B$ for some $a,b\in\RR$. 
If $V$ is also positive, and the spectral lower bound of $V|_{L^2(B)}$ is zero, then $V = a WE_B$, $a\geq 0$. 
\ethm
\proof
The first statement is immediate from our discussion and the relation
 \ben\label{ES}
- (f , W  f)_B + (f , M f)_B = - (f , L f)_B + ||f||^2_B\,  ,
\een
cf. \eqref{WL},
taking into account that the Friedrichs extensions of $-E_B L$ on $\S$ is equal to $-L_B E_B$ by Prop. \ref{Fri}.

$WE_B$ commutes with $\F_B$ by Prop. \ref{PropW}. 
The characterization of $V$ follows by an argument similar to the one in the proof of Prop. \ref{WF}. 
\endproof
The parabolic distribution $(1 - r^2)$ appears in both the parabolic and the Legendre entropy expression. 
Near the center of $B$, the parabolic entropy is close to the Born entropy. On the other hand, near the boundary of $B$, the prolate entropy gets close to the Born entropy. 

Let us specialize now on the one-dimensional case as studied in \cite{SP1} (on the even functions subspace of $L^2(B)$). As $T_B$ is strictly positive and Hilbert-Schmidt, its eigenvalues can be ordered as $\l_1 > \l_2 > \cdots  > 0$; moreover, they are simple; the eigenvalues of $-W_B$ can be ordered as $\a_1 < \a_2 < \cdots <\infty$; they correspond to the $\l_k$'s in inverse order, that is $T_B$ and $-W_B$ share the same $k$-the eigenfunction $f_k$, which is unique up to a phase once we normalize it as $||f_k||_B^2 = 1$. Then
\[
(f_k, T_B f_k)_B = \l_k\, ,\quad  - (f_k, W_B f_k)_B = \a_k\, ,
\]
and $\pi \a_k$ is the prolate entropy of $f_k$. 

As the information is the opposite of the entropy, the above relations show the intuitive fact that the functions with lower prolate entropy, thus higher information in $B$, are the ones with better support concentration in $B$ in space and Fourier modes. $f_1$ carries the best information as it is optimally concentrated. 

We expect the ordering correspondence between the eigenvalues of $T_B$ and $W_B$ to hold
 in the higher dimensional case too.  

\section{Concentration in balls of arbitrary radius}
We briefly indicate here how the results in this paper easily extend to the case of localization in balls of any radius. 
The more general prolate operator 
\[
W_{\rm min}(c) = \nabla (1 -r ^2) \nabla - c^2 r^2
\]
is studied in \cite{SP1}, $c > 0$. We consider $W_{\rm min}(c)$ as an operator on $L^2(\RR^d)$ with domain $\S$. 
Denote by $\d_\l$, $\l>0$, the dilation operator on $L^2(\RR^d)$
\[
(\d_\l f )(x) = \l^{-d/2} f(\l^{-1} x)\, ,
\]
so $\d_\l$ is a unitary operator. We also set $\F_\l = \d^{-1}_\l \F $; in particular, $\F_{2\pi}$ is the commonly used Fourier transform
in Communication Theory and elsewhere. 
\bprop
$W_{\rm min}(c)$ commutes with $\F_c$. 
\eprop
\proof
Since $\d^{-1}_\l r \d_\l = \l r$ and $\d^{-1}_\l \nabla \d_\l = \l^{-1} \nabla$, we have
\begin{align*}
 \F_c W_{\rm min}(c) \F_c^{-1} &=  \d^{-1}_c\F\big(  \nabla (1 -r ^2) \nabla - c^2 r^2 \big) \F^{-1} \d_c \\
&=  \d^{-1}_c \F\big(  \nabla (1 -r ^2) \nabla - r^2 + r^2 - c^2 r^2 \big) \F^{-1} \d_c \\
&=  \d^{-1}_c\F\big(  \nabla (1 -r ^2) \nabla - r^2\big)\F^{-1} \d_c +  \d^{-1}_c\F\big(r^2 - c^2 r^2 \big) \F^{-1} \d_c \\
&=  \d^{-1}_c\big(  \nabla (1 -r ^2) \nabla - r^2\big)\d_c  -  \d^{-1}_c \big(\nabla^2 -  c^2 \nabla^2 \big) \d_c \\
&=  c^{-2} \big(  \nabla (1 - c^{2} r ^2) \nabla \big) - c^{2} r^2 - c^{-2} \nabla^2 + \nabla^2  \\
&=  \nabla (c^{-2} - r ^2) \nabla - c^{2} r^2 - c^{-2} \nabla^2 +  \nabla^2  \\
&=   \nabla (1 - r ^2) \nabla  - c^2  r^2 = W_{\rm min}(c) \, .
\end{align*}
\endproof
The analysis of $W_{\rm min}(c)$ is now the same as in the case $c=1$. $W_{\rm min}(c)$ admits a natural extension $W_c$ that commutes with $\F_c$, $E_B$, $\hat E_{B_c}$, where $B_c$ denotes the ball of radius $c$ centered at the origin. 

The prolate operator corresponding to the localization in balls $B_\l$, and $B_{\l'}$ in Fourier transform, is obtained
by conjugating $W_c$ by the dilation operator, that is $W_{\l,\l'} = \d_\l W_c \d_\l^{-1}$, 
\[
W_{\l,\l'} = \nabla(\l^2 - r^2)\nabla -{ \l'}^2 r^2\, ,
\]
$\l\l' = c$.   
$W_{\l,\l'}$ commutes with 
$E_{B_\l} = \d_\l E_B \d_\l^{-1}$ and $\d_\l \hat E_{B_c} \d_\l^{-1} 
= \hat E_{B_{\l'}}$. 

Now, 
\[
W_{\l,\l'} 
= \nabla(\l^2 - r^2)\nabla  + \l^{-2} c^2 (\l^2 - r^2) - c^2 
= \nabla(\l^2 - r^2)\nabla  + {\l'}^2 M_\l - c^2
\]
on $L^2(B_\l)$, with $L_\l$ is the natural extension of $\nabla(\l^2 - r^2)\nabla$ and $M_\l = (1 - r^2)$, thus
\ben\label{genW}
- \pi (f, W_{\l,\l'} f)_{B_\l} + \l^2 {\l'}^2 \pi ||f||_{B_\l}^2 = -\pi (f, L_\l f)_{B_\l} +  {\l'}^2 \pi (f, M_\l f)_{B_\l}  \, ,
\een
an entropy relation that generalizes Thm. \ref{thme}. 

The first term on the left of \eqref{genW} is the prolate entropy of $f$ w.r.t. $B_\l$ and $B_{\l'}$. 
Note however that the Legendre entropy $-\pi(f, L_\l f)_{B_\l}$ does not depend on $\l'$; as $\l'\to0$, the prolate entropy approaches the Legendre entropy. 
\bigskip

\noindent
{\bf Acknowledgements.} 
We thank A. Connes and G. Morsella for valuable comments.  
\smallskip

\noindent
The author is supported by 
GNAMPA-INdAM. 
He acknowledges the Excellence Project 2023-2027 MatMod@TOV awarded to the Department of Mathematics, University of Rome Tor Vergata.

\end{document}